\documentclass[11pt]{amsart}
\usepackage{amsmath}
\usepackage{amsthm}
\usepackage{amsbsy}
\usepackage{amsopn}
\usepackage{amsmath,amscd}
\usepackage{amssymb}
\usepackage{amsfonts}
\usepackage{multirow}
\usepackage{mathrsfs} 
\usepackage[official ]{eurosym}
\usepackage[style=alphabetic, backend=bibtex, maxbibnames=99]{biblatex}
\newcommand{\RNum}[1]{\uppercase\expandafter{\romannumeral #1\relax}}

\newcommand\cF{\mathcal{F}}

\newcommand\cH{\mathcal{H}}
\newcommand\cI{\mathcal{I}}

\newcommand\cW{\mathcal{W}}

\def\bC{\mathbf{C}}

\def\bP{\mathbf{P}}

\def\bR{\mathbf{R}}

\def\bZ{\mathbf{Z}}

\newcommand{\GL}{\textup{GL}}
\newcommand{\Pic}{\textup{Pic}}

 \usepackage[all]{xy}
                \usepackage{amsthm,amsfonts,amsmath,amscd,amssymb,epsfig,verbatim,
}

\usepackage[parfill]{parskip}    
\usepackage{graphicx}
\usepackage{epstopdf}
\date{\today} 
 
\usepackage{graphicx}

  \defbibheading{bibliography}[\refname]{\section*{#1}}

\usepackage[dvipsnames,svgnames,x11names,hyperref]{xcolor}
\usepackage{url,graphicx,verbatim,amssymb,enumerate,stmaryrd}
\usepackage[colorlinks,citecolor=Bittersweet,linkcolor=Bittersweet,urlcolor=Bittersweet,filecolor=Bittersweet]{hyperref}
\usepackage{tikz} 
\usetikzlibrary{arrows,decorations.pathmorphing,backgrounds,fit,positioning,shapes.symbols,chains,calc}
\usepackage[all]{xy}
\xyoption{matrix}
\xyoption{arrow}
\usepackage[hmargin=3cm,vmargin=3cm]{geometry}
\usepackage{setspace,kantlipsum}
\tikzset{help lines/.style={step=#1cm,very thin, color=gray},
help lines/.default=.5} 
\usepackage{tikz-cd}
\usepackage{booktabs}
 
\theoremstyle{plain}
\newtheorem{theorem}{Theorem}[section]
\newtheorem{thm}[theorem]{Theorem}

\newtheorem{prop}[theorem]{Proposition}

\theoremstyle{definition}

\newtheorem{inductive step}[theorem]{Inductive step}

\newtheorem{inductive lemma}[theorem]{Inductive Lemma}
\theoremstyle{remark}

\newtheorem*{remark*}{Remark}

\newtheorem{question}[theorem]{Question}

\newtheorem*{example*}{Example}

\newcommand{\into}{\hookrightarrow}

\newcommand{\bS}{{\bf S}}

\renewcommand{\Re}{{\rm Re\,}}

\numberwithin{figure}{section}

\usepackage{stackengine}
\def\rlwd{.4pt}
\def\rlht{1.1pt}
\def\rhatvrule{\rule{\rlwd}{\rlht}}
\def\rhat#1{%
  \setbox0=\hbox{$#1$}%
  \stackon[0pt]{\stackon[1pt]{\ensuremath{#1}}{%
    \rhatvrule\kern\wd0\kern-\rlwd\kern-\rlwd\rhatvrule}}%
    {\rule{\wd0}{\rlwd}}%
}

\usepackage{tikz}
\usetikzlibrary{calc, arrows.meta}

\usepackage{mathtools}

\newcommand{\Sq}{\mathrm{Sq}}
\theoremstyle{plain}

\theoremstyle{remark}

\newif\ifnotes
\notestrue
\definecolor{notecolour}{rgb}{0.10,0.25,0.65}

\title[A Weinstein manifold without a spectral refinement]{A Weinstein manifold without any spherical Fukaya category}
\author{Abigail Ward}
\address{University of Sheffield, Sheffield, UK}
\email{abigail.ward1@sheffield.ac.uk}
\date{\today}

\begin{document}

\begin{abstract}
We show that, after attaching a subcritical handle to $T^*\bP^2$ to obtain a manifold $M$, the integral wrapped Floer endomorphism algebra of a cotangent fibre admits no spherical refinement for any choice of brane structure; this obstructs the existence of any wrapped Fukaya category for $M$ over the sphere spectrum. We further discuss some features of the resulting Weinstein eight-manifold $M$, which is homotopy equivalent to $S^4$.
\end{abstract}

\maketitle

\section{Introduction}

One central question in Floer homotopy theory asks how and when  the Fukaya category of a Weinstein manifold $X$ can be defined over the sphere spectrum $\bS$. Existing Floer-theoretic constructions use tangential data \cite{Lar21, PS24, ADP24, Bla24, PS25, PS25a, BB25} to provide the necessary coherent orientations on moduli spaces of pseudoholomorphic curves. The general expectation, borne out in work in progress of Asplund--Blakey--Deshmukh--Porcelli \cite{ABDP}, is that such a spherical Fukaya category can be defined in the presence of \emph{Maslov data}, i.e., a nullhomotopy of a certain classifying map
\[ m_X: X \to \, B^2\,\Pic\,\bS. \]
After constructing such a category over $\bS$, one obtains an integral wrapped Fukaya category by base change.

Cohomological obstructions to the existence of Maslov data were exhibited in \cite{AGLW}. These limit the reach of the available constructions, but do not in themselves rule out the existence of a spherical refinement. The following theorem does.

 \begin{thm} \label{thm: main} There exists a Weinstein eight-manifold $M \simeq S^4$ containing a cylindrical Lagrangian $L\cong \bR^4$ such that for any choice of brane structure on $L$, the Floer complex $CW^*(L,L)$ admits no spherical refinement: there is no unital $E_1$-ring spectrum
$A$ together with an equivalence $
H\bZ\otimes_{\bS} A\simeq CW^*(L,L;\bZ) $
of $E_1$-algebras over $H\bZ$. The Lagrangian $L$ generates the integral wrapped Fukaya category $\cW(M)$. 
 \end{thm}
 
This gives an example of the phenomenon sought in \cite[Question~1.7]{AGLW}: there exists a space in which the failure of Maslov data is detected by the Floer theory of $L$ itself, and  no `abstract' spherical Fukaya category containing $L$ as an object can exist.  To our knowledge, this is the first example of this sort. 
 
 The manifold $M$ is constructed by attaching a subcritical handle to $T^*\bP^2$ along an embedded $S^2$ representing the hyperplane class $[\bP^1]$.  The Lagrangian $L$ is a cotangent fibre chosen disjoint from the attaching region. Attaching a subcritical handle does not change the Floer complex of $L$, so from the cotangent-fibre computation of Abouzaid \cite[Thm.~1.1 and Rem.~1.2]{AbouzaidBasedLoops}, we compute that  
\[ CW_{M; b=0}^*(L,L) \simeq CW_{T^*\bP^2; b=0}^*(L, L) \simeq C_{-*}(\Omega \bP^2; \kappa),\]
where $\kappa$ is the nontrivial rank-one integral local system. The strategy of the proof is to show that this complex admits no compatible Steenrod operations.

 \subsection{Brane structures and background classes} The Maslov map  $m_X$ is defined as the composition
\[ X \xrightarrow{[TX]}  BU \xrightarrow{\beta} B^2(\bZ \times BU) \xrightarrow{B^3 J} B^2\,\Pic\,\bS, \]
where $\beta$ is the complex Bott periodicity map, $B^3 J$ is the thrice delooped complex $J$-homomorphism, and $\bS$ is the sphere spectrum.
More generally, for a commutative ring spectrum $R$, one says that $X$ admits $R$-Maslov data if the composition 
\[ X\xrightarrow{m_X}B^2\,\Pic\,\bS\to B^2\,\Pic\,R\]
 is nullhomotopic \cite[Sec.~4.1]{AGLW}, where the second map  is induced by extension of scalars along the unit map $\bS\to R$. One expects that a choice of $R$-Maslov data yields an $R$-linear Fukaya category. Evidently any choice of spherical Maslov data for $X$ will yield $R$-Maslov data, and one expects the construction in progress of Asplund--Blakey--Deshmukh--Porcelli to be compatible with this choice of Maslov data via base change $\cW(X;\bS) \to \cW(X;\bZ)$. 

After choosing spherical Maslov data $\mathfrak m$, one also needs brane data on the Lagrangians $L\subset X$ to construct spectral versions of Lagrangian Floer complexes. For a given $L$, this consists of a nullhomotopy of the map $L \to B \, \Pic\, \bS$ obtained by comparing $\mathfrak{m}|_{L}$ and the nullhomotopy of $m_X|_{L}$ furnished by the Lagrangian tangent bundle $TL$. 

In the integral case, choices of $H\bZ$-Maslov data recover the usual grading and orientation data for wrapped Fukaya categories, namely a trivialisation of $(\det_{\bC}TX)^{\otimes2}$ and a background class $b\in H^2(X;\bZ/2)$, while the brane data for $L$ consist of a grading and a relative Pin structure with background class $b$ \cite[Lem.~2.26]{ADPCocores}.  Our example makes the role of the choice of background data and brane structures particularly visible. Since $T^* \bP^2$ has a polarisation, it admits a spherical Fukaya category. The apparent contradiction is resolved by remembering that there are \emph{two} choices of background class on $T^*\bP^2$, since $H^2(T^*\bP^2;\bZ/2)\cong\bZ/2$. The background class $b=w_2(\bP^2)$ is the one obtained from the standard cotangent polarisation. Indeed, 
\[ CW_{b}^*(\mathrm{fibre},\mathrm{fibre}) \simeq C_{-*}(\Omega \bP^2),\]
which has a spherical refinement $\Sigma^\infty_+\Omega\bP^2$, as expected. 
After attaching the handle, only the $b=0$ choice survives. 

\subsection*{Related work}
Blakey--Porcelli \cite{BlakeyPorcelli} show that for odd $n\geq3$, the canonical $\mathrm{MU}$-valued symplectic cohomology of $T^*\bP^n$ admits no spherical lift; our construction is a multiplicative open-string analogue of this result. Blakey \cite{BlakeyViterbo} constructs a Thom-spectrum model for the $\mathrm{MU}$-valued symplectic cohomology of $T^*Q$ for any closed smooth $Q$ using an $\mathrm{MU}$-local system on $\mathcal{L}Q$. For $Q=\bP^2$, the integral orientation twist obtained by base changing from $\mathrm{MU}$ to $H\bZ$ restricts to $\kappa$ on $\Omega\bP^2$; the expected spectral open--closed and closed--open maps (see \cite{PS25a} for work on the former) should relate the corresponding open- and closed-string theories.

A related use of handle attachment appears in Husin \cite[Thm.~1.1 and Sec.~3]{Husin}, wherein one sees two closed exact Lagrangians whose Maslov classes cannot simultaneously vanish for any choice of grading data.

\subsection{Tangential data} 
In \cite{AGLW} one examines the relationship between (i) the geometry of singular Lagrangian skeleta, and in particular the possibility of finding an \emph{arboreal} skeleton for any given manifold, (ii) the topological properties of the tangent bundle, and (iii) the obstructions for Floer theory.  As with the obstructions for Maslov data, the obstructions to the existence of an arboreal skeleton are cohomological rather than geometrically effective; no explicit skeleton that cannot be arborealised is exhibited in that work. 

 This leads us to ask: 
\begin{question} \label{question: arb} Does $M$ admit an arboreal skeleton? \end{question}
A necessary condition for a Weinstein manifold $X$ to admit such a skeleton is the existence of an $\cI$-structure \cite[Thm.~3.34]{AGLW}. For $\dim_{\bR}X=2n$, this amounts to a lift $f$ in the diagram below:  
\[
\begin{tikzcd}
& \cI(n)
    \arrow[d] \\
X \arrow[r, "{[TX]}"'] \arrow[ur, dashed, "f"]
& BU(n)
\end{tikzcd}
\]
where the classifying space $\cI(n)$ is the homotopy pushout
\[ \cI(n) = BO(n)\underset{BO(n-1)\times BO(1)}{\cup^{h}}\bigl(BO(n-1)\times BU(1)\bigr).\]
$M$ admits an $\cI$-structure, so a negative answer to Question~\ref{question: arb} would disprove the naive $h$-principle for $\cI$-structures.  In any case, the results of \cite{AGEN} obstruct the existence of any positive skeleton. 

We study the tangential structures on $M$ in Sec.~3, and include a direct computation of its Maslov obstruction as an element of $\pi_4(B^2\,\Pic\,\bS)$. We also exhibit low-dimensional criteria for various structures, using $M$ to complete a list of examples showing the relevant dimension bounds are sharp.  

\subsection*{AI declaration} The example of $M$ and $L$ was found in conversation with GPT-6 Astra after the author asked for an example of a Floer complex without a spectral refinement and suggested possible ambient spaces; the LLM pointed out that the cotangent fibre in $T^* \bP^2$ already gave such an example. ChatGPT and Claude were also used to proofread the document.  

\subsection*{Acknowledgements} The author's understanding of this subject was shaped by many conversations with Daniel \'{A}lvarez-Gavela and Tim Large, and she is very grateful to both of them. She also thanks Kenneth Blakey and Thomas Kragh for conversations about orientations for $T^*\bP^2$; Noah Porcelli and Ivan Smith for helpful discussions about the statements in this note; and Kenneth Blakey for his thoughtful comments on a draft. The observation that every Weinstein manifold of dimension at most ten admits an $\cI$-structure is due to S{\o}ren Galatius, in conversations with John Pardon and the authors of \cite{AGEN}. This  work was completed during a visit to the Simons Center for Geometry and Physics for the programme “Interactions between quantum and homotopical geometry”; she thanks the SCGP for its hospitality.
The author was partially supported by UKRI Frontier Research Grant EP/X030660/1 held by Ivan Smith (in lieu of an ERC Advanced Grant) and by a Royal Society University Research Fellowship held by Cheuk Yu Mak.

\section{Construction of $M$ and $L$} We first construct our manifold $M$. Begin with $T^* \bP^2$, containing its unit cosphere bundle $S^* \bP^2$, and choose a smooth embedding 
\[ i: S^2 \into S^*\bP^2\]
of a sphere representing a generator of $\pi_2(T^* \bP^2)$.  By the $h$-principle for subcritical isotropic embeddings, we can isotope $i$ to a smooth isotropic embedding along which we can attach a Weinstein $3$-handle. The resulting manifold $M$ is homotopy equivalent to $\bP^2/\bP^1 \simeq S^4$. It appears also in work of Torricelli \cite[Ex.~3.4]{Torricelli}. 

Let $L \subset M$ be a cotangent fibre of $T^*\bP^2$, chosen disjoint from the attaching region of the handle. Consider now the wrapped complex $CW^*(L,L)$ where $L$ is taken to have the trivial spin structure. 

\begin{prop}\label{prop:HW}
We have 
\[
HW^{j}(L,L;\bZ)\cong
\begin{cases}
\bZ/2 & \text{if } j\in\{0,-4,-8,\dots\},\\
0 & \text{otherwise.}
\end{cases}
\]
\end{prop}
\begin{proof}
Since $M$ is obtained from $T^*\bP^2$ by attaching a subcritical handle away from $L$, \cite[Cor.~1.29]{GPS} gives
\[
HW^*(L,L;\bZ)\cong HW^*_{T^*\bP^2;\,b=0}(L,L;\bZ).
\]
Here the background class is zero, since $H^2(M;\bZ/2)=0$. By the cotangent-fibre computation of Abouzaid \cite[Thm.~1.1 and Rem.~1.2]{AbouzaidBasedLoops}, we have
\[
HW^*_{T^*\bP^2;\,b=0}(L,L;\bZ)\cong H_{-*}(\Omega\bP^2;\kappa),
\]
where $\kappa$ is the rank-one integral local system obtained by transgressing $w_2(\bP^2)$; see also work of Seidel \cite[Ex.~1 and Rem.~3]{SeidelCP2} for a discussion of this orientation twist, following Kragh. The local system $\kappa$ has monodromy $-1$ around a generator of $\pi_1(\Omega\bP^2) \cong\bZ$. 

The Hopf fibration gives a homotopy equivalence
\[
\Omega\bP^2\simeq S^1\times\Omega S^5;
\]
see \cite[Sec.~3]{CohenJonesYan}. Under this equivalence $\kappa$ is pulled back from the nontrivial local system $\rho$ on $S^1$, so by the K\"unneth formula we have
\[ H_*(\Omega\bP^2;\kappa)\cong H_*(S^1;\rho)\otimes_{\bZ}H_*(\Omega S^5;\bZ). \]
The homology $H_*(S^1;\rho)$ is concentrated in degree $0$, where it is given by $\bZ/2$. 
By \cite[Sec.~3, eq.~(3.3)]{CohenJonesYan}, we have
\[ H_*(\Omega S^5;\bZ)\cong\bZ[u],\qquad |u|=4. \]
The statement in the proposition follows. 
\end{proof}

One can also compute the ring structure, but it is not needed for our purposes. We can now prove our main theorem, that $CW^*(L,L)$ has no spherical refinement.
\begin{proof}[Proof of Thm.~\ref{thm: main}]Since $M\simeq S^4$, and in particular $H^1(M;\bZ)=H^2(M;\bZ/2)=0$, it admits a unique choice of integral Maslov data; since $L$ is contractible, its brane structure is unique up to a shift in grading. It therefore suffices to rule out a spherical refinement of the algebra computed in Prop.~\ref{prop:HW}.
We do this using the classical obstruction to a unital multiplication on the mod-two Moore spectrum \cite{Barratt, Rognes}.  

Suppose to a contradiction that there exists an $E_1$-ring spectrum $A$ with
$A \otimes_{\bS} H\bZ \simeq CW^*(L,L;\bZ)$. Then  $A$ carries a unital multiplication $A \wedge A \to A$ and has $H_*(A;\bZ) = \pi_*(A \otimes_{\bS} H\bZ) \cong HW^{-*}(L,L;\bZ)$. Since $A$ is a spectrum, we know that $H_*(A;\bZ/2)$ carries natural Steenrod operations
\[ \Sq_*^i:H_*(A;\bZ/2)\longrightarrow H_{*-i}(A;\bZ/2), \]
which satisfy the Cartan formula with respect to the product induced by the multiplication on $A$.  Calculating with the universal coefficient theorem and the above proposition gives
\[ H_j(A;\bZ/2)\cong
\begin{cases}
\bZ/2 & \text{if } j\geq 0 \text{ and } j\equiv 0,1\pmod{4},\\
0 & \text{otherwise.}
\end{cases}
\]

Let $x\in H_1(A;\bZ/2)$ be the nonzero class. The integral Bockstein map gives an isomorphism  
\[ H_1(A;\bZ/2) \cong H_0(A;\bZ)\cong\bZ/2,\]
 sending $x$ to the unit; reducing modulo two gives $\Sq^1_*(x)=1$. Note that $x^2$ and $\Sq_*^2(x)$ vanish for degree reasons. The Cartan formula gives 
\[ 0=\Sq_*^2(x^2) =\Sq_*^2(x)x+\Sq_*^1(x)\Sq_*^1(x)+x\Sq_*^2(x) =1, \]
a contradiction. \end{proof}

\section{Tangential structures on $M$}
The calculations above show that the Maslov obstruction  \[ [m_M] \in \pi_4(B^2\,\Pic\,\bS) \cong \bZ/2\]  is necessarily nonvanishing, as otherwise a choice of $\bS$-Maslov data would induce the unique choice of $H\bZ$-Maslov data on $M$ and yield a spherical refinement of $CF^*(L,L)$ by \cite{ADPCocores}, contradicting the preceding calculation. One can also compute the class $[m_M]$ directly. The universal Maslov map induces the following maps on homotopy groups:
\[
\begin{tikzcd}[column sep=large, row sep=large]
\pi_4(BU)
  \arrow[r] \arrow[d, "\cong"']
&
\pi_4(B(U/O))
  \arrow[r] \arrow[d, "\cong"']
&
\pi_4(B^2\,\Pic\,\bS)
  \arrow[d, "\cong"]
\\
\bZ
  \arrow[r]
&
\bZ/2
  \arrow[r, "\cong"']
&
\bZ/2.
\end{tikzcd}
\]

The second arrow is the stable $J$-homomorphism
$\pi_1(O)\to\pi_1^s$, which is an isomorphism. Here we identify $\pi_4(BU)$ with $\bZ$ by the Chern number 
$ [f]\mapsto \langle f^*c_2,[S^4]\rangle.$

Let $i:\bP^2\hookrightarrow M$ denote the original zero-section, and let $h\in H^2(\bP^2;\bZ)$ be the hyperplane class. The inclusion induces an isomorphism $i^*:H^4(M;\bZ)\to H^4(\bP^2;\bZ)$. We have
\[ i^*c_2(TM) = c_2(T_{\bR}\bP^2\otimes\bC) = -p_1(T_{\bR}\bP^2) = -3h^2,  \]
 so  the parity of this class shows that $[m_M] \in \bZ/2$ is nonvanishing as expected.

But $M$ has an $\cI$-structure; in fact, we have the following low-dimensional results: 

\begin{prop}\label{prop:low-dimensional}
Let $X$ be a Weinstein manifold of real dimension $2n$.
\begin{enumerate}
\item[(a)] If $2n\le10$, then $X$ admits an $\cI$-structure.
\item[(b)] If $2n\le10$, then $X$ admits a polarisation if and
only if it admits spherical Maslov data.
\item[(c)] If $2n\le6$, then $X$ admits a polarisation if and
only if $2c_1(X)=0$. \end{enumerate} 
\end{prop}

\begin{proof} (a) Fix an $\omega$-compatible almost complex structure $J$ and set $g=\omega(\cdot,J\cdot)$.  Consider a family of symmetric complex bilinear forms on the tangent bundle of $X$, i.e., a section $q$ of the bundle $\mathrm{Sym}_{\bC}^2(T^*X)$. We claim that $q$ can be chosen to have corank at most one at every point. Indeed, the locus of forms of corank at least two has complex codimension three in each fibre; recall that in general the codimension of corank $k$ forms is $k(k+1)/2$. Thus the complement is $4$-connected; since $X$ has the homotopy type of a CW complex of dimension at most five, the standard obstruction theory argument shows that $q$ exists. 

Define a real self-adjoint linear operator $B$ on the underlying real bundle $TX$ by $\Re q(v,w) = g(Bv,w)$. We now define a dual $\cI$-structure as in \cite{AGLW}, following an idea in \cite[proof of Prop.~4.2, p.~21]{Seidel}.  At each point $x \in X$ declare $\cF_x^+ \subset T_x X$ to be the sum of the positive eigenspaces of $B$. On the locus $U$ where $q$ has full rank, the resulting distribution $\cF^+$ is Lagrangian. On $V_0 = X \setminus U$, define a rank $n+1$ distribution $\cH $ by taking the $g$-orthogonal complement of $\cF^+$; this is automatically transverse to $\cF^+$ and co-isotropic. We can then take a small neighbourhood $V$ containing $V_0$ and extend $\cH$ as a co-isotropic distribution. After possibly shrinking $V$, this extension is transverse to $\cF^+$ on $U \cap V$, yielding a dual $\cI$-structure subordinate to the cover $(U,V)$. 

(b) One only has to prove the reverse direction. Factor the Maslov map of $X$ as the composition of the two maps \[f: X \to BU \to B(U/O); \quad \quad \quad \phi: B(U/O) \to B^2\,\Pic\,\bS.\] By \cite[Sec.~4.1]{AGLW}, $\phi$ fits into a homotopy commutative diagram
\[
\begin{tikzcd}[column sep=large, row sep=large]
B^3 O \arrow[r] \arrow[d, "B^3 J"']
  & B(U/O) \arrow[r] \arrow[d, "\phi"]
  & B^2\bZ \arrow[d, equal] \\
B^3 \GL_1(\bS) \arrow[r]
  & B^2\,\Pic\,\bS \arrow[r]
  & B^2\bZ,
\end{tikzcd}
\]
whose rows are fibrations. Since $B^2\bZ=K(\bZ,2)$ and both fibres are $2$-connected, the associated long exact sequences on homotopy groups show that both $B(U/O)$ and $B^2\,\Pic\,\bS$ are simply connected and that the projections to $B^2\bZ$ induce isomorphisms on $\pi_2$; thus $\phi$ induces an isomorphism on $\pi_2$ as well.

For $k\geq 3$, the groups $\pi_k(B^2\bZ)$ and $\pi_{k+1}(B^2\bZ)$ vanish, so the inclusions of the fibres induce isomorphisms
\[
\pi_k(B(U/O))\cong\pi_k(B^3O)\cong\pi_{k-3}(O), \qquad
\pi_k(B^2\,\Pic\,\bS)\cong\pi_k(B^3\GL_1(\bS))\cong\pi_{k-3}(\GL_1(\bS)).
\]
Under these identifications, the map induced by $\phi$ is the $J$-homomorphism
\[
J_*:\pi_{k-3}(O)\longrightarrow\pi_{k-3}(\GL_1(\bS)).
\]
On $\pi_0$, this sends a component of $O$ to its determinant in $\{\pm1\}$, and is therefore an isomorphism. On $\pi_1$, it sends the generator of $\pi_1(O)\cong\bZ/2$ to the nonzero stable Hopf class $\eta\in\pi_1^s\cong\bZ/2$, so it is again an isomorphism; see \cite[Ex.~6.20]{Ranicki}. Thus $\phi$ induces isomorphisms on homotopy groups in degrees at most four. Finally,
\[
\pi_5(B(U/O))\cong\pi_2(O)=0.
\]

Let $F$ be the homotopy fibre of $\phi$. The calculation above shows that $\pi_k(F) = 0$ for $k \leq 3$, and the inclusion $i$ induces the zero map on $\pi_4$. When $\phi \circ f$ is nullhomotopic, we can write $f = i \circ g$ for some lift $g: X \to F$. We now use the fact that $X$ has the homotopy type of a CW complex of dimension at most five. We can choose $g$ to be constant on the $3$-skeleton $X_3 \subset X$; the obstruction to extending the nullhomotopy of $i \circ g$ across each $4$-cell lies in the image of $\pi_4(F)$ under $i_*$, which is trivial; the obstruction to extending across each $5$-cell lies in $\pi_5(B(U/O))$, which is trivial as well. This gives a stable polarisation, which destabilises to a polarisation of $TX$ because $X$ has the homotopy type of an $n$-dimensional CW complex. 

(c) Complex vector bundles on spaces of dimension at most $3$ are stably classified by $c_1$, so $2c_1(X)=0$ implies that $TX$ has a stable polarisation; but then a stable polarisation implies the existence of a polarisation by the usual dimensional argument, e.g. \cite[Sec.~1.1]{AGLW}.
\end{proof}

These bounds are sharp: the twelve-dimensional example in \cite[Prop.~4.9 and Rem.~4.11]{AGLW} admits spherical Maslov data but no $\cI$-structure, and the manifold $M$ constructed above has $c_1 =0$ but no polarisation.

\printbibliography

@thesis{Lar21,
  author       = {Large, Tim},
  title        = {Spectral {Fukaya} Categories for {Liouville} Manifolds},
  type         = {phdthesis},
  institution  = {Massachusetts Institute of Technology},
  year         = {2021},
  url          = {https://hdl.handle.net/1721.1/139233}
}

@article{PS24,
  author       = {Porcelli, Noah and Smith, Ivan},
  title        = {Spectral {Floer} theory and tangential structures},
  journaltitle = {Transactions of the American Mathematical Society},
  volume       = {379},
  year         = {2026},
  doi          = {10.1090/tran/9594},
  eprint       = {2411.03257},
  eprinttype   = {arxiv},
  eprintclass  = {math.SG}
}

@online{ADP24,
  author       = {Asplund, Johan and Deshmukh, Yash and Pieloch, Alex},
  title        = {Spectral equivalence of nearby {Lagrangians}},
  year         = {2024},
  eprint       = {2411.08841},
  eprinttype   = {arxiv},
  eprintclass  = {math.SG}
}

@online{Bla24,
  author       = {Blakey, Kenneth},
  title        = {{Floer} homotopy theory and degenerate {Lagrangian} intersections},
  year         = {2024},
  note         = {To appear in Journal of Symplectic Geometry},
  eprint       = {2410.11478},
  eprinttype   = {arxiv},
  eprintclass  = {math.SG}
}

@online{PS25,
  author       = {Porcelli, Noah and Smith, Ivan},
  title        = {Bordism from quasi-isomorphism},
  year         = {2025},
  eprint       = {2509.21587},
  eprinttype   = {arxiv},
  eprintclass  = {math.SG}
}

@online{PS25a,
  author       = {Porcelli, Noah and Smith, Ivan},
  title        = {Open-closed maps and spectral local systems},
  year         = {2025},
  eprint       = {2509.21483},
  eprinttype   = {arxiv},
  eprintclass  = {math.SG}
}

@online{BB25,
  author       = {Blakey, Kenneth and Bonciocat, Ciprian Mircea},
  title        = {Parameterized {Lagrangian} {Floer} homotopy},
  year         = {2025},
  eprint       = {2506.20122},
  eprinttype   = {arxiv},
  eprintclass  = {math.SG}
}

@online{BlakeyViterbo,
  author       = {Blakey, Kenneth},
  title        = {Spectral {Viterbo} isomorphism: complex-oriented versus framed},
  year         = {2026},
  eprint       = {2608.20289},
  eprinttype   = {arxiv},
  eprintclass  = {math.SG}
}

@online{BlakeyPorcelli,
  author       = {Blakey, Kenneth and Porcelli, Noah},
  title        = {Bulk-deformations, {Floer} complex bordism, and {Grothendieck--Riemann--Roch}},
  year         = {2026},
  eprint       = {2605.06620},
  eprinttype   = {arxiv},
  eprintclass  = {math.SG}
}

@online{Husin,
  author       = {Husin, Axel},
  title        = {{Maslov} class of exact {Lagrangians} and cylindrical handles},
  year         = {2025},
  eprint       = {2502.14750},
  eprinttype   = {arxiv},
  eprintclass  = {math.SG}
}

@online{ADPCocores,
  author       = {Asplund, Johan and Deshmukh, Yash and Pieloch, Alex},
  title        = {Homotopy rigidity of nearby {Lagrangian} cocores},
  year         = {2025},
  eprint       = {2511.09548},
  eprinttype   = {arxiv},
  eprintclass  = {math.SG}
}

@unpublished{ABDP,
  author       = {Asplund, Johan and Blakey, Kenneth and Deshmukh, Yash and Porcelli, Noah},
  title        = {Spectral wrapped {Fukaya} categories},
  note         = {In preparation},
  shorthand    = {ABDP}
}

@unpublished{AGLW,
  author       = {{\'A}lvarez-Gavela, Daniel and Large, Tim and Ward, Abigail},
  label        = {AGLW},
  title        = {On arborealization, {Maslov} data, and lack thereof},
  year         = {2025},
  note         = {To appear in Journal of Topology},
  eprint       = {2503.09783},
  eprinttype   = {arxiv},
  eprintclass  = {math.SG},
  url          = {https://arxiv.org/abs/2503.09783}
}

@unpublished{AGEN,
  author       = {{\'A}lvarez-Gavela, Daniel and Eliashberg, Yakov and Nadler, David},
  label        = {AGEN},
  title        = {Positive arborealization of polarized {Weinstein} manifolds},
  year         = {2020},
  note         = {Preprint},
  eprint       = {2011.08962},
  eprinttype   = {arxiv},
  eprintclass  = {math.SG},
  url          = {https://arxiv.org/abs/2011.08962}
}

@article{GPS,
  author       = {Ganatra, Sheel and Pardon, John and Shende, Vivek},
  title        = {Sectorial descent for wrapped {Fukaya} categories},
  journaltitle = {Journal of the American Mathematical Society},
  year         = {2024},
  volume       = {37},
  number       = {2},
  pages        = {499--635},
  doi          = {10.1090/jams/1035},
  eprint       = {1809.03427},
  eprinttype   = {arxiv},
  eprintclass  = {math.SG}
}

@article{AbouzaidBasedLoops,
  author       = {Abouzaid, Mohammed},
  title        = {On the wrapped {Fukaya} category and based loops},
  journaltitle = {Journal of Symplectic Geometry},
  year         = {2012},
  volume       = {10},
  number       = {1},
  pages        = {27--79},
  doi          = {10.4310/JSG.2012.v10.n1.a3},
  eprint       = {0907.5606},
  eprinttype   = {arxiv},
  eprintclass  = {math.SG}
}

@unpublished{SeidelCP2,
  author       = {Seidel, Paul},
  title        = {A remark on the symplectic cohomology of cotangent bundles, after {Kragh}},
  date         = {2010-09},
  note         = {Informal note},
  url          = {https://math.mit.edu/~seidel/texts/cp2.pdf}
}

@incollection{CohenJonesYan,
  author       = {Cohen, Ralph L. and Jones, John D. S. and Yan, Jun},
  title        = {The loop homology algebra of spheres and projective spaces},
  booktitle    = {Categorical Decomposition Techniques in Algebraic Topology},
  editor       = {Arone, Gregory and Hubbuck, John and Levi, Ran and Weiss, Michael},
  series       = {Progress in Mathematics},
  volume       = {215},
  pages        = {77--92},
  publisher    = {Birkh{\"a}user},
  location     = {Basel},
  year         = {2003},
  doi          = {10.1007/978-3-0348-7863-0_5},
  eprint       = {math/0210353},
  eprinttype   = {arxiv},
  eprintclass  = {math.AT}
}

@article{Torricelli,
  author       = {Torricelli, Brunella Charlotte},
  title        = {Projective twists and the {Hopf} correspondence},
  journaltitle = {Algebraic \& Geometric Topology},
  year         = {2024},
  volume       = {24},
  number       = {8},
  pages        = {4139--4200},
  doi          = {10.2140/agt.2024.24.4139},
  eprint       = {2006.12170},
  eprinttype   = {arxiv},
  eprintclass  = {math.SG}
}

@thesis{Seidel,
  author       = {Seidel, Paul},
  title        = {{Floer} homology and the symplectic isotopy problem},
  type         = {D.Phil. thesis},
  institution  = {University of Oxford},
  year         = {1997},
  url          = {https://math.mit.edu/~seidel/texts/thesis.pdf}
}

@article{Barratt,
  author = {Barratt, M. G.},
  title = {Track groups ({II})},
  journaltitle = {Proceedings of the London Mathematical Society},
  series = {3},
  volume = {5},
  number = {3},
  date = {1955},
  pages = {285--329},
  doi = {10.1112/plms/s3-5.3.285}
}

@article{Rognes,
  author = {Rognes, John},
  title = {The product on topological {Hochschild} homology of the integers with mod 4 coefficients},
  journaltitle = {Journal of Pure and Applied Algebra},
  volume = {134},
  number = {3},
  date = {1999},
  pages = {211--218},
  doi = {10.1016/S0022-4049(97)00154-0}
}

@book{Ranicki,
  author = {Ranicki, Andrew},
  title = {Algebraic and Geometric Surgery},
  series = {Oxford Mathematical Monographs},
  publisher = {Oxford University Press},
  location = {Oxford},
  date = {2002},
  isbn = {9780198509240},
  doi = {10.1093/acprof:oso/9780198509240.001.0001}
}

\end{document}